\documentclass[aps,pre,amsmath,amsfonts,amssymb,superscriptaddress,showpacs,11pt]{revtex4}
\usepackage{graphics}
\usepackage[latin1]{inputenc}
\usepackage{amscd}
\usepackage{soul}

\usepackage[pdftex]{graphicx}

\usepackage{rotating}
\usepackage[pdftex,hypertexnames=true,linktocpage=true]{hyperref}

\usepackage{color}

\usepackage{amsthm}
\usepackage{framed}
\usepackage{amssymb}
\usepackage{amsmath,mathrsfs}
\usepackage{hhline}
\usepackage{ulem}
\usepackage{fancyhdr}
\usepackage{tikz}
\everymath{\displaystyle}
\everymath{\displaystyle}

\everymath{\displaystyle}
\newcommand{\beq}{\begin{equation}}
\newcommand{\eeq}{\end{equation}}

\newcommand{\bi}{\begin{itemize}}
\newcommand{\ei}{\end{itemize}}
\newcommand{\D}{\mathrm{D}}

\renewcommand{\H}{\mathcal{H}}

\newcommand{\K}{\mathrm{K}}

\def\RR{\mathbb{R}}

\newcommand{\Tau}{\mathcal{T}}

\newcommand{\Ma}{\mathrm{M}}
\newcommand{\esse}{\mathcal{S}}

\newcommand{\Div}{\mathcal{D}}
\newcommand{\nihat}{\mathrm{X}}
\newcommand{\grad}{\mathrm{grad}}
\newcommand{\total}{\mathrm{d}}
\newcommand{\RC}{\mathrm{R}}

\newcommand{\Riemann}{\mathcal{R}}
\newcommand{\tangent}{{\rm T}}

\newcommand{\paralleltransport}{{\rm P}}

\newcommand{\metric}{{\rm g}}
\newcommand{\btheta}{\boldsymbol{ \theta}}
\newcommand{\bxi}{\boldsymbol{ \xi}}

\newcommand{\bieta}{\boldsymbol{ \eta}}

\begin{document}

\title{\textbf{Divergence functions in Information Geometry}}
\author{Domenico Felice}
\email{domenico.felice@unicam.it}
\affiliation{Max Planck Institute for Mathematics in the Sciences\\
 Inselstrasse 22--04103 Leipzig,
 Germany}
\author{Nihat Ay}
\email{nay@mis.mpg.de}
\affiliation{ Max Planck Institute for Mathematics in the Sciences\\
 Inselstrasse 22--04103 Leipzig,
 Germany\\
Santa Fe Institute, Santa Fe, NM 87501, USA\\
 Faculty of Mathematics and Computer Science, University of Leipzig, PF 100920, 04009 Leipzig, Germany}

\begin{abstract}
A recently introduced canonical divergence $\Div$ for a dual structure $(\metric,\nabla,\nabla^*)$ on a smooth manifold $\Ma$ is discussed in connection to other divergence functions. Finally, open problems concerning symmetry properties are outlined.
\end{abstract}


\maketitle
\section{Introduction}

The geometrical structure induced by a {\it divergence function} (or {\it contrast function}) on a smooth manifold $\Ma$ provides a unified 
approach to measurement of notions as information, energy, entropy, playing an important role in mathematical sciences to research random phenomena \cite{Eguchi92}. In the mathematical formulation, a divergence function $\Div(p,q)$ on a smooth manifold $\Ma$ is defined by the first requirement for {a} distance:
\begin{equation}
\label{distance}
\Div(p,q)\geq 0\,,\qquad \Div(p,q)=0\quad\mbox{iff}\quad p=q\ .
\end{equation}
An important example of a divergence function is given by the Kullback-Leibler divergence $\K(p,q)$ in the context that $p$ and $q$ are the vectors of probabilities of disjoint events \cite{Eguchi85}, namely
\begin{equation}
\label{K-Ldivergence}
\K(p,q)=\sum_{i=1}^{n+1}\,p_i\,\log\left(\frac{p_i}{q_i}\right)
\end{equation}
is a function on the $n$-simplex $\esse:=\{p=(p_1,\ldots,p_{n+1})\,|\,p_i\geq 0\,,\sum_i p_i=1\}$. Given a smooth $n$-dimensional manifold $\Ma$, we assume that $\Div:\Ma\times\Ma\rightarrow\RR^+$ is a $C^{\infty}$-differentiable function. Working with the local coordinates $\{\bxi_p:=(\xi_p^1,\ldots,\xi^n_p)\}$ and $\{\bxi_q:=(\xi_q^1,\ldots,\xi^n_q)\}$ at $p$ and $q$, respectively, it follows from Eq. (\ref{distance}) that
\begin{align}
\label{firstderivative}
& \left.\partial_i\,\Div(\bxi_p,\bxi_q)\right|_{p=q}=0,\quad \left.\partial^{\prime}_i\,\Div(\bxi_p,\bxi_q)\right|_{p=q}=0\\
\label{secondderivative}
& \left.\partial_j\partial_i\,\Div(\bxi_p,\bxi_q)\right|_{p=q}=-\left.\partial^{\prime}_j\partial_i\,\Div(\bxi_p,\bxi_q)\right|_{p=q}=\left.\partial^{\prime}_j\partial^{\prime}_i\,\Div(\bxi_p,\bxi_q)\right|_{p=q}\,,
\end{align}
where $\partial_i=\frac{\partial}{\partial\xi^i_p}$ and $\partial_i^{\prime}=\frac{\partial}{\partial\xi^i_q}$. Moreover, under the assumption that
\begin{equation}
\label{positivity}
\metric_{ij}:=\left.\partial_i\partial_j\Div(\bxi_p,\bxi_q)\right|_{p=q}>0
\end{equation}
we can see that the manifold $\Ma$ is endowed, through the divergence function $\Div$, with the Riemannian metric tensor given by $\metric=\metric_{ij}\,\total\xi^i\otimes\total\xi^j$, where the Einstein notation is adopted. The symmetry of $\metric$ immediately follows from the requirement that $\Div$ is a $C^{\infty}$ function on $\Ma\times\Ma$.

From Eq. (\ref{K-Ldivergence}) we can see that, in general, a divergence function $\Div$ is not symmetric.
The asymmetry of $\Div$ leads to two different affine connections, $\nabla$ and $\nabla^*$, on $\Ma$ such that $1/2(\nabla+\nabla^*)$ is the Levi-Civita connection with respect to the metric tensor $\metric=\metric_{ij}\total\xi^i\otimes\total\xi^j$ defined by Eq. \eqref{positivity}. More precisely, working with the local coordinates $\{\bxi_p\}$ and $\{\bxi_q\}$, we can define the symbols $\Gamma_{ijk}$ and $\Gamma_{ijk}^*$ of the connections $\nabla$ and $\nabla^*$, i.e. $\Gamma_{ijk}=\metric\left(\nabla_{\partial_i}\partial_j,\partial_k\right)$ and ${\Gamma}^*_{ijk}=\metric\left(\nabla^*_{\partial_i}\partial_j,\partial_k\right)$,  by means of the following relations
\begin{equation}
\label{connections}
\Gamma_{ijk}(p)=-\left.\partial_i\partial_j\partial_k^{\prime} \Div(\bxi_p,\bxi_q)\right|_{p=q}, \qquad {\Gamma}^*_{ijk}(p)=-\left.\partial^{\prime}_i\partial^{\prime}_j\partial_k \Div(\bxi_p,\bxi_q)\right|_{p=q}\ .
\end{equation} 
To sum up, a divergence function $\Div$ on a smooth manifold $\Ma$ induces a metric tensor on $\Ma$ by Eq. \eqref{positivity}. In addition, the divergence $\Div$ yields two linear torsion-free connections, $\nabla$ and $\nabla^*$, on $\tangent\Ma$ which are dual with respect to the metric tensor $\metric$ \cite{Eguchi85}:
\begin{equation}
\label{dualstructure}
X\,\metric\left(Y,Z\right)=\metric\left(\nabla_X Y,Z\right)+\metric\left(Y,\nabla^*_X Z\right),\ \forall\ X,Y,Z\in\Tau(\Ma),
\end{equation}
where $\Tau(\Ma)$ denotes the space of vector fields on $\Ma$. Finally, we refer to the quadruple $(\Ma,\metric,\nabla,\nabla^*)$ as a {\it statistical manifold} \cite{Ay17}.

\subsection{The inverse problem within Information Geometry}

\noindent The inverse problem is to find a divergence $\Div$ which generates a given geometrical structure $(\Ma,\metric,\nabla,\nabla^*)$. For any such statistical manifold there exists a divergence $\Div$ such that Eq. (\ref{positivity}) and Eq. (\ref{connections}) hold true \cite{Matumoto93}. However, this divergence is not unique and there are infinitely many divergences generating the same geometrical structure $(\metric,\nabla,\nabla^*)$. When this structure is dually flat, namely the curvature tensors of $\nabla$ and $\nabla^*$ are null ($\RC(\nabla)\equiv 0$ and $\RC^*(\nabla^*)\equiv 0$), Amari and Nagaoka introduced a canonical divergence which is a Bregman divergence \cite{Amari00}. The canonical divergence has nice properties such as the
generalized Pythagorean theorem and the geodesic projection theorem and it turns out to be of uppermost importance to define a canonical divergence in the general case. A first attempt to answer this fundamental issue is provided by Ay and Amari in \cite{Ay15} where a canonical divergence for a general statistical manifold $(\Ma,\metric,\nabla,\nabla^*)$ is given by using the geodesic integration of the inverse exponential map. This one is understood as a {\it difference vector} that translates $q$ to $p$ for all $q,p\in\Ma$ sufficiently close to each other.

\begin{figure}
\centering
\begin{tikzpicture}
\draw [very thick] (0,0) ellipse (2 and 1);
\coordinate [label=above left:$\mathbb{R}^n$] (t) at (1.5,.2);
\draw [<-, very thick] (-.95,-.45) -- (0.3,.5)
node at (-0.8,-.65) {$p$}
node at (0.5,.5) {$q$}
node [pos=.5, sloped, above] {$p-q$};
\draw [black, thick] plot [mark=*, only marks]
coordinates {(-1.,-.5) (.3,.5)};
\node at (0,-1.6) {$(\textbf{A})$};
\end{tikzpicture}
\qquad
\begin{tikzpicture}
\draw [very thick] (0,0) ellipse (2 and 1.2);
\coordinate [label=above left:$\Ma$] (t) at (1.5,.5);
\draw [<-, very thick] (0,-.65) -- (0.3,.5);
\node (p) at (-0.8,-.7) {$p$};
\node (q) at (0.5,.5) {$q$};
\coordinate [label=below right:{\small $\nihat(q,p)$}] (X) at (0.1,0);
\draw [-latex, bend left, very thick] (0.3,.5) to (-1,-.5)
node at (-.4,-.05) {$\widetilde{\sigma}$};
\draw [black, thick] plot [mark=*, only marks]
coordinates {(-1.,-.5) (.3,.5)};
\draw [dashed] (-2,0) to [out=40, in=140] (2,0);
\node at (0,-1.6) {$(\textbf{B})$};
\end{tikzpicture}
\\
\vspace{.2cm}
\begin{tikzpicture}
\draw [very thick] (0,0) ellipse (2 and 1.2);
\coordinate [label=above left:$\Ma$] (t) at (1.5,.5);
\draw [->, very thick] (-1,-.5) -- (0,-.5);
\draw [->, very thick] (0.3,.5) -- (1.2,.4);
\node (p) at (-1.2,-.6) {$p$};
\node (q) at (0.3,.7) {$q$};
\coordinate [label=below right:{\small $\nihat(p,q)$}] (X) at (-1,-.6);
\draw [-latex, bend left, very thick] (-1,-.5) to (0.3,.5) 
node at (-.4,.45) {$\sigma^*$};
\coordinate [label=below right:{\small $\Pi_q(p)$}] (X) at (.3,.4);
\draw [-latex, bend right, very thick] (-1,-.5) to (0.3,.5) 
node at (-.2,-.05) {$\sigma$};
\draw [black, thick] plot [mark=*, only marks]
coordinates {(-1.,-.5) (.3,.5)};
\draw [dashed] (-2,0) to [out=50, in=130] (2,0);
\node at (0,-1.6) {$(\textbf{C})$};
\end{tikzpicture}
\caption{Illustration $(\mathbf{A})$ of the difference vector in $\RR^n$ pointing from $q$ to $p$; and $(\mathbf{B})$ the difference vector $\nihat(q,p)=\dot{\widetilde{\sigma}}(0)$ as the inverse of the exponential map in $q$. The novel vector $\Pi_q(p)$ at $q$ is described in $(\mathbf{C})$ as the $\nabla$ parallel transport of $\nihat(p,q)=\dot{\sigma}(0)$ along the $\nabla^*$-geodesic $\sigma^*$ from $p$ to $q$.}\label{vectors}
\end{figure}
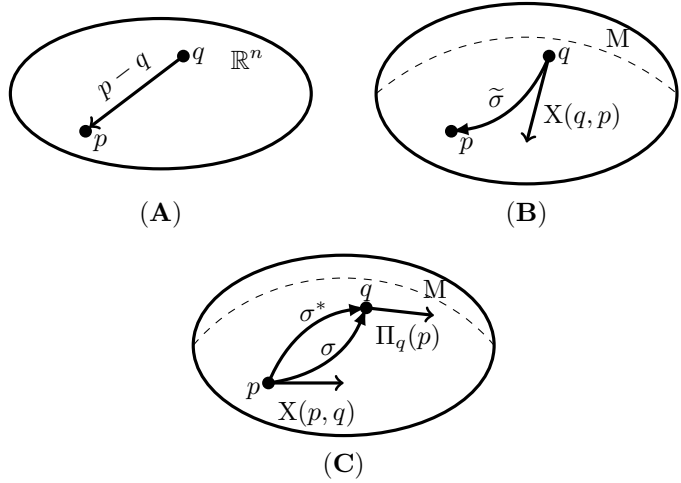

To be more precise, the inverse exponential map supplies a generalization to $\Ma$ of the concept of difference vector in $\RR^n$. In detail, let $p,q\in\RR^n$, the difference between $p$ and $q$ is the vector $p-q$ pointing to $p$ (see side ($\mathbf{A}$) of Fig. \ref{vectors}). Then, the difference between $p$ and $q$ in $\Ma$ is provided by the inverse exponential map. In particular, given $p,q$ suitably close in $\Ma$, the difference vector from $q$ to $p$ is defined as (see ($\mathbf{B}$) of Fig. \ref{vectors})
\begin{equation}
\label{Nihatvector}
\nihat_q(p):=\nihat(q,p):=\exp_q^{-1}(p)=\dot{\widetilde{\sigma}}(0)\,,
\end{equation}
where $\widetilde{\sigma}$ is the $\nabla$-geodesic from $q$ to $p$. Therefore, the divergence $\D$ introduced in \cite{Ay15} is defined as the path integral
\begin{equation}
\label{AyDiv}
\D(p,q):=\int_0^1\, \langle \nihat_t(p),\dot{\widetilde{\sigma}}(t)\rangle_{\widetilde{\sigma}(t)}\ dt\,,\quad \nihat_t(p):=\nihat(\widetilde{\sigma}(t),p):=\exp_{\widetilde{\sigma}(t)}^{-1}(p)\,,
\end{equation}
where  $\langle\cdot,\cdot\rangle_{\widetilde{\sigma}(t)}$ denotes the inner product induced by $\metric$ on $\widetilde{\sigma}(t)$. After elementary computation Eq. \eqref{AyDiv} reduces to,
\begin{equation}
\label{Aydivergence}
\D(p,q)=\int_0^1 \, t\, \|\dot{\sigma}(t)\|^2\ dt\ ,
\end{equation}
where $\sigma(t)$ is the $\nabla$-geodesic from $p$ to $q$ \cite{Ay15}. The divergence $\D(p,q)$ has nice properties such as the positivity and it reduces to the canonical divergence proposed by Amari and Nagaoka when the manifold $\Ma$ is dually flat. However, if we consider definition \eqref{AyDiv} for a general path $\gamma$, then $\D_{\gamma}(p,q)$ will be depending on $\gamma$. On the contrary, if the vector field $\nihat_t(p)$ is integrable, then $\D_{\gamma}(p,q)=: \D(p,q)$ turns out to be independent of the path from $q$ to $p$. 

\section{A new canonical divergence}

\noindent In this article, we discuss about a novel divergence function recently introduced in \cite{Felice18}. This turns out to be a generalization of the divergence introduced by Ay and Amari. The definition of the new divergence {(see below Eqs. \eqref{canonicaldivergence}, \eqref{dualcanonicaldivergence})} relies {on an extended analysis of the intrinsic structure} of the dual geometry of a general statistical manifold $(\Ma,\metric,\nabla,\nabla^*)$. In particular, we introduced a vector at $q$ {by modifying} the definition (\ref{Nihatvector}) of the difference vector $\nihat(q,p)$. Consider $p,q\in\Ma$ such that there exist both, a unique $\nabla$-geodesic $\sigma$ and a unique $\nabla^*$-geodesic $\sigma^*$, connecting $p$ to $q$. Moreover, let $\nihat_p(q):=\exp_p^{-1}(q)= \dot{\sigma}(0)\in\tangent_p\Ma$, we then $\nabla$-parallel translate it along $\sigma^*$ from $p$ to $q$ (see ($\mathbf{C}$) of Fig. \ref{vectors}), and obtain 
\begin{equation}
\label{P}
{\Pi_q(p):=}\paralleltransport_{\sigma^*}\nihat_p(q)\in\tangent_q\Ma\ .
\end{equation} 
(Note that $\Pi_q(p)$ corresponds to minus a difference vector.) Analogously, we introduce the dual vector of $\Pi_q(p)$ as the $\nabla^*$-parallel transport of $\dot{\sigma}^*(0)$ along the $\nabla$-geodesic $\sigma$,
\begin{equation}
\label{P*}
\Pi_q^*(p):=\paralleltransport^*_{\sigma}\nihat_p^*(q),\quad \nihat_p^*(q):=\stackrel{*}{\exp}_p^{-1}(q) = \dot{\sigma}^*(0)\ ,
\end{equation}
where $\stackrel{*}{\exp}_p$ denotes the exponential map of the $\nabla^*$-connection. A fundamental result obtained in \cite{Felice18} is that the sum of $\Pi_q(p)$ and $\Pi^*_q(p)$ is the gradient of a symmetric function $r_p(q)=r(p,q)$ that we call {\it pseudo-norm}:
\begin{equation}
\label{MainResult}
\Pi_q(p)+\Pi^*_q(p)=\grad_q r_p,\quad r(p,q):=\langle\exp_p^{-1}(q),\stackrel{*}{\exp}_p^{-1}(q)\rangle_p\,.
\end{equation}
By letting $q\in\Ma$ be varying, we then obtain two vector fields whenever $p$ and $q$ are connected by a unique $\nabla$-geodesic and a unique $\nabla^*$-geodesic. Then, we can introduce two vector fields on an arbitrary path $\gamma:[0,1]\rightarrow\Ma$ connecting $p$ and $q$ in the following way. Let us firstly assume that for each $t\in[0,1]$ there exist a unique $\nabla$-geodesic $\sigma_t$ and a unique  $\nabla^*$-geodesic $\sigma^*_t$ connecting $p$ {with} $\gamma(t)$. Then, we define
\begin{align}
\Pi_t(p)= \paralleltransport_{\sigma^*_t} \nihat_p(t),&\qquad \nihat_p(t)=\exp_p^{-1}(\gamma(t))\label{vectorfieldPi}\\
\Pi_t^*(p)=\paralleltransport^*_{\sigma_t}\nihat_p^*(t),&\qquad \nihat_p^*(t)=\stackrel{*}{\exp}_p^{-1}(\gamma(t))\ . \label{vectorfieldPi*}
\end{align} 
Therefore, from Eq. \eqref{MainResult} we have that the sum
\begin{equation}
\label{indipendence}
\int_0^1\ \langle \Pi_t(p),\dot{\gamma}(t)\rangle_{\gamma(t)}\ dt+\int_0^1\  \langle \Pi_t^*(p),\dot{\gamma}(t)\rangle_{\gamma(t)}\ dt= r_p(q)
\end{equation}
is independent of the particular path from $p$ to $q$.

At this point, we define the novel canonical divergence $\Div$, and its dual function $\Div^*$, from $p$ to $q$ by the geodesic integration of $\Pi_t(p)$ and $\Pi^*_t(p)$, respectively. In particular, we have that
\begin{align}
\label{canonicaldivergence}
& \Div(p,q):=\Div_p(q):=\int_0^1\,\langle\Pi_t(p),\dot{\sigma}(t)\rangle_{\sigma(t)}\,\total t\,,\\
\label{dualcanonicaldivergence}
& \Div^*(p,q):=\Div^*_p(q):=\int_0^1\,\langle\Pi^*_t(p),\dot{\sigma}^*(t)\rangle_{\sigma^*(t)}\,\total t\,,
\end{align}
where $\sigma$ and $\sigma^*$ are the $\nabla$-geodesic and the $\nabla^*$-geodesic  from $p$ to $q$, respectively.

In this manuscript we review {the relation} of the canonical divergence $\Div(p,q)$ {to other divergence functions} in Section \ref{SectionII}. Finally, we outline in Section \ref{SectionIII} the open problems concerning the symmetry properties of $\Div$.

\section{Comparison with previous divergence functions}\label{SectionII}

Given a general statistical manifold $(\Ma,\metric,\nabla,\nabla^*)$, the {basic requirement} for a smooth function $\Div:\Ma\times\Ma\rightarrow\RR$ to be a divergence on $\Ma$ is its consistency with the dual structure $(\metric,\nabla,\nabla^*)$ through Eqs. (\ref{positivity})-(\ref{connections}) and the positivity $\Div(p,q)>0$ for all $p,q\in\Ma$ sufficiently close to each other such that $p\neq q$. The novel canonical divergence (\ref{canonicaldivergence}) succeeds to holding these properties (see Theorem 5 in \cite{Felice18}).

\begin{figure}[h!]
\centering
\begin{tikzpicture}
\draw  (0,0) to [out=60, in=120] (5,3);
\coordinate [label=below left:$p$]
(p) at (0,0);
\coordinate [label=right:$q$] (q) at (5,3);
\coordinate [label=above left:$\sigma(t)$] (t) at (1.8,2.2);
\draw (0,0) to [out=-30,in=-60] (2,2.55)
       node at (2.2,.8) {$\sigma^*_t$}
node at (0.5,1.3) {$\sigma_t$} ;
\draw [black, thick] plot [mark=*, only marks]
coordinates {(0,0) (5,3) (2,2.55)};
\draw [->, ultra thick] (2,2.55) -- (1.7,3.2);
\coordinate [label=below left:{\small $\Pi_t(p)$}] (P) at (2,3.6); 
\draw [->, ultra thick] (0,0) -- (.5,.8);
\draw [->, ultra thick] (0,0) -- (.9,-.5);
\coordinate [label=left:{\small $\nihat_p(t)$}] (P) at (0.3,0.6); 
\draw [->, ultra thick] (2,2.55) -- (2.9,3.3);
\coordinate [label=below left:{\small $\nihat^*_p(t)$}] (P) at (1,-.4);
\coordinate [label=below left:{\small $-\nihat_t(p)$}] (P) at (3.5,3.0);
\end{tikzpicture}
\caption{The vector $\Pi_t(p)$ is obtained by $\nabla$-parallel translating the vector $\nihat_p(t):=\dot{\sigma}_t(0)$ along the $\nabla^*$-geodesic $\sigma^*_t$. The opposite of the difference vector $\nihat_t(p)$ at $\sigma(t)$ can be understood as the $\nabla$-parallel translation of the vector $\nihat_p(t)$ along the $\nabla$-geodesic $\sigma_t$.}\label{Pi&X}
\end{figure}
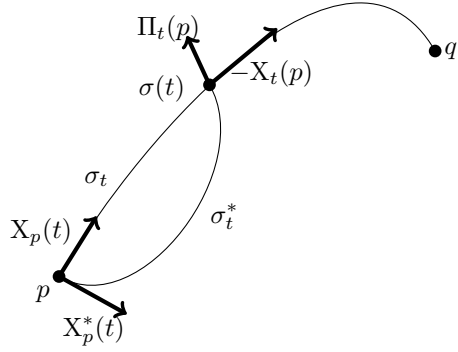

In this section, {we will show that the canonical divergence $\Div$ can be interpreted as a generalization of the divergence $\D$ introduced by Ay and Amari. Indeed, we will see that these two divergences coincide on particular classes of statistical manifolds. In order to achieve this result, we investigate some geometric properties of the vector field $\Pi_t(p)$ given by Eq. (\ref{vectorfieldPi}) aiming to split such a vector field in terms of the difference vector $\nihat_t(p)$ given in Eq. (\ref{AyDiv}). To be more precise, let us refer to Fig. \ref{Pi&X} where the $\nabla$-geodesic $\sigma(t)\,(0\leq t\leq 1)$ connecting $\sigma(0)=p$ with $\sigma(1)=q$ is drawn. Then,} 
for each $t\in[0,1]$ we can consider the $\nabla$-geodesic $\sigma_t(s)\,(0\leq s\leq 1)$ connecting $p$ {with} $\sigma(t)$ and the $\nabla^*$-geodesic $\sigma^*_t(s)\,(0\leq s\leq 1)$ connecting $p$ {with} $\sigma(t)$. The difference vector $\nihat_t(p)=\nihat(\sigma(t),p)$ at $\sigma(t)$ pointing to $p$ is given in terms of the inverse exponential map by
$
\nihat_t(p):=\exp^{-1}_{\sigma(t)}(p).
$
Therefore, the opposite of $\nihat_t(p)$ can be viewed as the $\nabla$-parallel translation of $\nihat_p(t)=\exp_p^{-1}(\sigma(t))$ along the $\nabla$-geodesic $\sigma_t$, namely 
$
-\nihat_t(p)=\paralleltransport_{\sigma_t}\nihat_p(t)\ .
$
Consider now the loop $\Sigma_t$ based at $p$ and given by first traveling from $p$ to $\sigma(t)$ along the $\nabla^*$-geodesic $\sigma_t^*$ and then back from $\sigma(t)$ to $p$ along the reverse of the $\nabla$-geodesic $\sigma_t$. If $\Sigma_t$ lies in a sufficiently small neighborhood of $p$, then \cite{Felice18}
\begin{equation*}
\label{Parallel&Curvature}
\paralleltransport_{\Sigma_t} \nihat_p(t)=\nihat_p(t)+ \Riemann_{\Sigma_t}\left(\nihat_p^*(t),\nihat_p(t)\right),
\end{equation*}
where 
\begin{equation*}
\label{RiemannSigmat}
\Riemann_{\Sigma_t}\left(\nihat_p^*(t),\nihat_p(t)\right) :=\int_{B_t} \frac{\paralleltransport\left[ \Riemann\left(\nihat^*(t),\nihat(t)\right)\nihat(t)\right]}{\|\nihat_p^*(t)\wedge\nihat_p(t)\|}\ \total A
\end{equation*}
with $\nihat^*(t)$ and $\nihat(t)$ being the $\nabla$-parallel transport of $\nihat_p^*(t)\,(:=\stackrel{*}{\exp_p}^{-1}(\sigma(t)))$ and $\nihat_p(t)$, respectively, from $p$ to each point of $B_t$ along the unique $\nabla$-geodesic joining them. 
Here, $\Riemann$ is the curvature tensor of $\nabla$, $B_t$ denotes the dis{c} defined by the curve $\Sigma_t$ and $\nihat^*_p(t)$, $\nihat_p(t)$ are linearly independent. In addition, $\paralleltransport$ within the integral denotes the $\nabla$-parallel translation from each point in $B_t$ to $p$ along the unique $\nabla$-geodesic segment joining them. Finally, by means of the property of the parallel transport, we obtain the following geometric relation between the vector $\Pi_t(p)$ and the opposite of the difference vector $\nihat_t(p)$ \cite{Felice18},
\begin{align}\label{Pi&nihat}
\Pi_t(p)&= \paralleltransport_{\sigma_t}\nihat_p(t)+\ \paralleltransport_{\sigma_t}\left[\Riemann_{\Sigma_t}\left(\nihat_p^*(t),\nihat_p(t)\right)\right]\ .
\end{align}
By noticing that 
$
\paralleltransport_{\sigma_t}\nihat_p(t)=\dot{\sigma}_t(1)=t\ \dot{\sigma}(t)
$
and inserting Eq. (\ref{Pi&nihat}) into the definition (\ref{canonicaldivergence}) of $\Div$, we obtain
\begin{equation}
\label{decompositiondivergence}
\Div(p,q)=\D(p,q)+\int_0^1 \ \langle\paralleltransport_{\sigma_t}\left[\Riemann_{\Sigma_t}\left(\nihat_p^*(t),\nihat_p(t)\right)\right],\dot{\sigma}(t)\rangle_{\sigma(t)}\ \total t\ ,
\end{equation}
where $\D(p,q)$ is the divergence introduced in \cite{Ay15} and given by Eq. (\ref{Aydivergence}). 

{It is clear from Eq. (\ref{decompositiondivergence}) that particular conditions on the curvature tensor would lead to the required equivalence between $\Div$ and $\D$. Actually, in Information Geometry classes of statistical manifolds are characterized by the conditions on the curvature tensors of $\nabla$ and $\nabla^*$ (see for instance Refs. \cite{Amari00}, \cite{Lauritzen87}, \cite{Zhang07}).} In the Table \ref{Tab1} we can see the categories of statistical manifolds on which the canonical divergence $\Div$ reduces to the divergence $\D$ introduced in \cite{Ay15}.  A statistical manifold $(\Ma,\metric,\nabla,\nabla^*)$ is self-dual when $\nabla=\nabla^*$. Therefore, in this case $\Ma$ becomes a Riemannian manifold endowed with the Levi-Civita connection. Hence, the vectors $\nihat_p(t)$ and $\nihat^*_p(t)$ coincide for all $t\in[0,1]$. Finally, the {skew-symmetry} of the curvature tensor $\Riemann$ yields the property $\Riemann(X,X)=0$ for any $X\in\Tau(\Ma)$. When a manifold $\Ma$ is dually flat, it has a mutually dual affine coordinates $\{\btheta\},\,\{\bieta\}$ and two potentials $\{\psi,\varphi\}$ such that $\D(p,q)=\psi(\btheta_p)-\varphi(\bieta_q)+\btheta_p\cdot\bieta_q$\cite{Ay15}. This claims that the canonical divergence $\Div$ coincides with the canonical divergence of Bregman type introduced in \cite{Amari00} by Amari and Nagaoka on dually flat manifolds. The concept of a symmetric statistical manifold, that is the information geometric analogue to {a} symmetric space in Riemannian geometry, was introduced in \cite{Henmi}. Here, the authors employed the following conditions on the curvature tensor, $\nabla\Riemann=0$ and $\RC(Y,X,X,X):=\langle\Riemann(Y,X)X,X\rangle=0$, in order to prove that their divergence function is independent of the particular path connecting any two points $p,q$ sufficiently close to each other. The connection between the canonical divergence $\Div$ and the divergence introduced by Henmi and Kobayashi is widely discussed in \cite{Felice18}. 

\begin{table}
\centering
\begin{tabular}{|ll|ll|}
\hline
\hskip 0.7truecm Statistical manifold \hskip 0.7truecm & & & \hskip 0.7truecm Condition on $\Riemann$ \hskip 0.7truecm \\
\hline
\hline
\hskip 0.7truecm Self-dual \hskip 0.7truecm & & & \hskip 0.7truecm $\Riemann(X,X)=0$ \hskip 0.7truecm \\ 
\hline
\hskip 0.7truecm Dually flat\hskip 0.7truecm & & &\hskip 0.7truecm $\Riemann\equiv 0$ \hskip 0.7truecm\\ 
\hline
\hskip 0.7truecm Symmetric \hskip 0.7truecm& & &\hskip 0.7truecm $\nabla\Riemann\equiv 0\, , \RC(Y,X,X,X)=0$  \hskip 0.7truecm\\ 
\hline
\end{tabular}
\vspace{.1cm}
\caption{The column on the left describes the category of statistical manifolds on which the canonical divergence $\Div$ reduces to the divergence $\D$. The column on the right shows the properties of the curvature tensor characterizing the corresponding manifolds and supplying the equivalence $\Div\equiv\D$ through the Eq. \eqref{decompositiondivergence}. $X,Y\in\Tau(\Ma)$.}\label{Tab1}
\end{table}

{To summarize, Tab. \ref{Tab1} describes, from the top to the bottom, the statistical manifolds ordered from less generality to more generality where the equivalence between $\Div$ and $\D$ is achieved. In this view, we can consider $\Div$ as an extension of the divergence $\D$ to the very general statistical manifold $(\Ma,\metric,\nabla,\nabla^*)$.}

Since Eq. (\ref{K-Ldivergence}) we know that in general a divergence function is not symmetric in its argument. However, the symmetry property owned by the canonical divergence of Bregman type on dually flat manifolds, namely $\D(q,p)=\D^*(p,q)$, shows the way for the further investigation about symmetry properties of $\Div$ in the very general context of Information Geometry.

\section{Future developments towards symmetry}\label{SectionIII} 

{The target of this section would be the description of the symmetry property $\Div(q,p)=\Div^*(p,q)$ for any statistical manifold $(\Ma,\metric,\nabla,\nabla^*)$. To this aim, we rely} on the {\it gradient--based approach} to divergence which was introduced in \cite{Ay15} and further developed in \cite{Ay17}. This approach yields  the  following decompositions of $\Pi_t(p)$ and $\Pi_t^*(p)$ in terms of the canonical divergence gradient and its dual \cite{Felice18},
\begin{align}
\label{Pigrad}
& \Pi_t(p)=\grad_{\gamma(t)}\Div_p(\gamma(t))+V_t,\quad \langle V_t,\dot{\sigma}_t(1)\rangle_{\gamma(t)}=0\\
\label{Pi*grad}
& \Pi^*_t(p)=\grad_{\gamma(t)}\Div^*_p(\gamma(t))+V^*_t,\quad \langle V^*_t,\dot{\sigma}^*_t(1)\rangle_{\gamma(t)}=0\,,
\end{align}
where $\sigma_t(s)\,, \sigma_t^*(s)\,,(0\leq s\leq 1)$ are the $\nabla$ and $\nabla^*$ geodesics, respectively, from $p$ to $\gamma(t)$ for any arbitrary path $\gamma(t)$ connecting $p$ and $q$.

\begin{figure}
\centering
\begin{tikzpicture}
\draw [very thick] (0,0) ellipse (2 and 1.2);
\coordinate [label=above left:$\Ma$] (t) at (1.5,.5);
\draw [very thick] (0.3,.5) to [out=270,in=65] (.1,-.8);
\node (H) at (-0.6,.2) {$\H$};
\node (q) at (0.5,.25) {$q$};
\coordinate [label=below right:{\small $L_q$}] (X) at (0.2,-.2);
\draw[very thick] (-1,.3) to [out=30, in=165] (.3,.5)
to [out=-10, in=130] (1.5,-.1);
\draw [black, thick] plot [mark=*, only marks]
coordinates {(.3,.5)};
\draw [dashed] (-2,0) to [out=40, in=140] (2,0);
\draw [thin] (.05,.55)-- (0.05,.25);
\draw [thin] (.05,.25)-- (0.3,.25);
\end{tikzpicture}
\caption{According to the Eguchi's theory \cite{Eguchi92}, the minimum contrast leaf $L_q$ turns out to be orthogonal at $q$ to the level-set $\H$ generated by the canonical divergence $\Div$.}\label{MinimumContrastLeaf}
\end{figure}
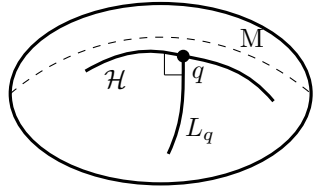

On the other hand, by means of the theory of minimum contrast geometry by Eguchi \cite{Eguchi92}, we can show that $\grad_{\gamma(t)}\Div_p$ is parallel to the tangent vector of the $\nabla$-geodesic starting from $p$ and $\grad_{\gamma(t)}\Div^*_p$ is parallel to the tangent vector of the $\nabla^*$-geodesic starting from $p$. This proves that Eqs. (\ref{Pigrad}) and (\ref{Pi*grad}) supply orthogonal decompositions of $\Pi_t(p)$ and $\Pi_t^*(p)$, respectively. To see this, let us consider the level sets of $\Div_p$ and $\Div^*_p$:
\begin{align}
\label{Hyper}
& \H(\kappa)=\{q\in\Ma\,|\,\Div_p(q)=\kappa\},\quad \H^*(\kappa)=\{q\in\Ma\,|\,\Div^*_p(q)=\kappa\}\,.
\end{align}
Then to each $q\in\H(\kappa)$ we can define the minimum contrast leaf of $\Div$ at $q$ \cite{Eguchi85}:
\begin{equation}
\label{Leaf}
L_q:=\{p\in\Ma\,|\,\Div(p,q)=\min_{q^{\prime}\in\H}\,\Div(p,q^{\prime})\}\,.
\end{equation}
Let us now fix $q$. Since $q$ minimizes the set $\{\Div(p,q)\,|\,p\in L_q\,,\,q\in\H\}$ it follows that the derivative of $\Div(p,q)$ at $q$ along any direction $U$ tangent to $\H$ vanishes, namely
$$
\partial^{\prime}_U\Div(p,q)=0,\quad \forall\, U\in\tangent_q\H\,,\, p\in L_q\,,
$$
where $\partial^{\prime}_U$ denotes the derivative at $q$ along the direction $U$. Thus we have that $\langle U,\Xi\rangle_q=0$ for all $U\in\tangent_q\H$ and $\Xi\in\tangent_q L_q$, or equivalently that the tangent space of $L_q$ coincides with the normal space of $\H$ at $q$ (see Fig. \ref{MinimumContrastLeaf} for a cross-reference). In addition, by taking derivatives at $q$ along directions $\Xi\,,\Upsilon$ normal to $\H$ we have that \cite{Eguchi92}
\begin{equation}
\label{II}
\metric\left(II(\Upsilon,\Xi),U\right):=-\left.\partial_{\Xi}\partial_{\Upsilon}\partial^{\prime}_U\Div(p,q)\right|_{p=q}=0\,,\quad \forall\, \Xi,\Upsilon\in\tangent_q L_q\,,\, U\in\tangent_q\H\,,
\end{equation}
where the first relation defines the second fundamental tensor $II$ with respect to the $\nabla$-connection. This implies that the second fundamental tensor with respect to $\nabla$ for $L_q$ vanishes at $q$. Therefore, according to the well-known Gauss formula \cite{Lee97}
$$
\nabla_{\Xi}\Upsilon=\nabla^{L_q}_{\Xi}\Upsilon+II(\Xi,\Upsilon)
$$
we can see from Eq. (\ref{II}) that the family of all curves which are orthogonal to the level set $\H$ are all $\nabla$-geodesics ending at $q$ (with a suitable choice of the parameter). 

Analogously, we have that the family of all curves which are orthogonal to the level set $\H^*$ are all $\nabla^*$-geodesics ending at $q$ (with a suitable choice of the parameter).

In order to answer what the relation of $\Div(p,q)$ and $\Div(q,p)$ is, let us put $\widetilde{\Div}(p,q):=\Div(q,p)$. Then, by noticing that $-\partial_i\partial_j\partial_k^{\prime}\,\widetilde{\Div}(p,q)=\Gamma^*_{ijk}$ and repeating the same arguments as above, we can show that $\grad_{q}\widetilde{\Div}_p=c(q)\,\grad_q\Div^*_p$, where $c:\Ma\rightarrow\RR$ is a smooth function on $\Ma$. This proves that there exists a function $f:[0,\infty]\rightarrow\RR^+$ and $f^{\prime}>0$ such that \cite{Felice18}
$$\Div(q,p)=f\left(\Div^*(p,q)\right)\ .$$ 
Though this relation holds for a very general statistical manifold $(\Ma,\metric,\nabla,\nabla^*)$, this result is still not satisfactory. {However, in Tab. \ref{Tab2} we can see the classes of statistical manifolds where the relation $\Div(q,p)=\Div^*(p,q)$ holds. This occurs in dually flat manifolds analogously to the canonical divergence of Bregman type introduced in \cite{Amari00}. Moreover, the required symmetry also holds in the symmetric statistical manifolds, which constitutes a new result in the setting of Information Geometry \cite{Felice18}. Forthcoming investigation will address such a symmetry in the general case.} 

\begin{table}
\centering
\begin{tabular}{|ll|ll|}
\hline
\hskip 0.7truecm Statistical manifold \hskip 0.7truecm & & & \hskip 0.7truecm Relation of $\Div(q,p)$ and $\Div^*(p,q)$ \hskip 0.7truecm \\
\hline
\hline
\hskip 0.7truecm Dually flat \hskip 0.7truecm & & & \hskip 0.7truecm $\Div(q,p)=\Div^*(p,q)$ \hskip 0.7truecm \\ 
\hline
\hskip 0.7truecm Symmetric\hskip 0.7truecm & & &\hskip 0.7truecm $\Div(q,p)=\Div^*(p,q)$ \hskip 0.7truecm\\ 
\hline
\hskip 0.7truecm General \hskip 0.7truecm& & &\hskip 0.7truecm $\Div(q,p)=f\left(\Div^*(p,q)\right)$  \hskip 0.7truecm\\ 
\hline
\end{tabular}
\vspace{.1cm}
\caption{The symmetry property of $\Div$ fails for general statistical manifolds. However, it holds true for dually flat and statistical symmetric manifolds \cite{Felice18}.}\label{Tab2}
\end{table}

%
%

\end{document}